\documentclass[12pt,leqno]{amsart}
\usepackage{amssymb,amsthm}
\usepackage{amsmath,amssymb,color}
\numberwithin{equation}{section}

\newcommand{\ep}{\varepsilon}
\newcommand{\la}{\lambda}
\newcommand{\va}{\varphi}
\newcommand{\ppp}{\partial}

\newcommand{\www}{\widetilde}

\newcommand{\M}{\mathbb{M}}
\newcommand{\N}{\mathbb{N}}

\newcommand{\ssss}{\sharp}
\newcommand{\ooo}{\overline}
\newcommand{\OOO}{\Omega}

\newcommand{\sumk}{\sum_{k=1}^{\infty}}

\newcommand{\sumij}{\sum_{i,j=1}^n}

\newcommand{\hhalf}{\frac{1}{2}}
\newcommand{\DDD}{\mathcal{D}}

\allowdisplaybreaks
\newcommand{\nuA}{\ppp_{\nu_A}}

\title
[]
{
Inverse parabolic problem with initial data by a single 
measurement
}
\author{
$^1$ O.Y. ~Imanuvilov and $^{2,3,4}$ M.~Yamamoto}

\thanks{
$^1$ Department of Mathematics, Colorado State University,\\
101 Weber Building, Fort Collins, CO 80523-1874, USA \\
e-mail: {\tt oleg@math.colostate.edu}
\\
$^2$ Graduate School of Mathematical Sciences, The University
of Tokyo, Komaba, Meguro, Tokyo 153-8914, Japan \\
$^3$ Honorary Member of Academy of Romanian Scientists,
Ilfov, nr. 3, Bucuresti, Romania \\
$^4$ Correspondence member of Accademia Peloritana dei Pericolanti,\\
Palazzo Universit\`a, Piazza S. Pugliatti 1 98122 Messina Italy \\
e-mail: {\tt myama@ms.u-tokyo.ac.jp}
}

\date{}
\begin{document}
\maketitle

\baselineskip 18pt

\begin{abstract}{\it 
We consider 
%an inverse coefficient problem for 
initial boundary value problems with the homogeneous Neumann 
boundary condition.  Given an initial value, we establish
the uniqueness in determining a spatially varying coefficient 
of zeroth-order term by a single measurement of Dirichlet data
on an arbitrarily chosen subboundary.  The uniqueness holds
in a subdomain where the initial value is positive, 
provided that it is sufficiently smooth which is specified 
by decay rates of the Fourier coefficients.
The key idea is the reduction to an inverse elliptic problem
and relies on elliptic Carleman estimates.}
\\
{\bf Key words.}
inverse coefficient problem, parabolic equation, uniqueness,
initial boundary value problem, inverse elliptic problem,
Carleman estimate
\\
{\bf AMS subject classifications.}35R30, 35K15
\end{abstract}

\section{Introduction}

Let $\Omega$ be a bounded domain in $\Bbb R^n$ with $C^2-$ boundary
$\ppp\OOO$. 
We set $\nu=(\nu_1(x),\dots,\nu_n(x))$ be an outward unit normal vector to 
$\ppp\OOO$, and
$$
\ppp_k = \frac{\ppp}{\ppp x_k}, \quad 1\le k\le n, \quad
\ppp_t = \frac{\ppp}{\ppp t}.
$$
We assume 
$$
a_{ij}=a_{ji} \in C^2(\ooo{\OOO}) \quad \mbox{for all $1\le i,j\le n$},
\quad c_1, c_2 \in C(\ooo{\OOO}),
$$
and there exists a constant $\kappa_1>0$ such that
$$
\sum_{i,j=1}^n a_{ij}(x)\xi_i\xi_j
\ge \kappa_1\vert \xi\vert^2 \quad \mbox{for all $x \in \ooo{\OOO}$ and
$\xi=(\xi_1, ... \xi_n) \in \Bbb R^n$}.
$$
We set 
$$
\ppp_{\nu_A}v: = \sumij a_{ij}(x)\nu_i(x)\ppp_jv(x), \quad 
x\in \ppp\OOO.
$$
Moreover let 
$$
A(x,D)v(x) = - \sumij \ppp_i(a_{ij}(x)\ppp_jv(x)) \quad 
\mbox{with} \quad
\DDD(A):= \{ v\in H^2(\OOO);\, \ppp_{\nu_{A_0}} v= 0\quad
\mbox{on $\ppp\OOO$}\}.                \eqno{(1.1)}
$$

In this article, we consider the folowing 
\\
{\bf Inverse problem.}
\\
{\it
Let 
$$\left\{ \begin{array}{rl}
& \ppp_tu = -A(x,D)u + c_1(x)u \quad \mbox{in $\OOO\times (0,T)$}, \\
& \nuA u = 0 \quad\mbox{on $\ppp\OOO$}, \\
& u(x,0) = a(x), \quad x\in \OOO,
\end{array}\right.
                             \eqno{(1.2)}
$$
and
$$\left\{ \begin{array}{rl}
& \ppp_tv = -A(x,D)v + c_2(x)v \quad \mbox{in $\OOO\times (0,T)$}, \\
& \nuA v = 0 \quad\mbox{on $\ppp\OOO$}, \\
& v(x,0) = a(x), \quad x\in \OOO.
\end{array}\right.
                                    \eqno{(1.3)}
$$
Let an initial value $a$ be suitably given and
$\gamma \subset \ppp\OOO$ be an arbitrarily chosen non-empty connected
relatively open subset of $\ppp\OOO$.
Then 
$$
\mbox{$u=v$ on $\gamma \times (0,T)$ implies $c_1=c_2$ in $\OOO$?} 
$$
}

This inverse problem has been intensively studied in the literature. 
Most general results are obtained for the case when a time of observation 
$t_0$ belongs to the  open interval $(0,T)$. In this case,  based on the method introduce in 
Bukhgeim and Klibanov \cite{BK}, Imanuvilov and Yamamoto \cite{IY98} proved 
the uniqueness and the Lipschitz stability in determination of coefficients 
corresponding to the zeroth order term.  
Recently in Imanuvilov and Yamamoto \cite{IY23},
the authors proved  a conditional Lipschitz stability 
estimate as well as the uniqueness for the case $t_0=T$.
See also Huang, Imanuvilov and Yamamoto \cite{HIY}.

In case an  observation is taken at the initial moment $t_0=0,$ 
to the authors' best knowledge, the question of uniqueness of a solution 
of inverse problem is open in general.
Limited to the one-dimensional case, a result in this direction was 
obtained by Suzuki \cite{S} and Suzuki and Murayama \cite{SM}. 
Klibanov \cite{Kl} proved the uniqueness of determination of zeroth order 
term coefficient in the case when $a_{ij}=\delta_{ij}$ (the case of 
the Laplace operator). Also it is assumed that the domain  of observation 
$\Gamma=\partial\Omega$.
The method proposed in \cite{Kl} is based on an integral transform and 
subsequent reduction of the original problem to the problem of determination 
of a coefficient of zeroth order term for a hyperbolic equation. 
After that, the method in \cite{BK} is applied.  
It should be mentioned that the method introduced by 
Bukhgeim and Klibanov is based on Carleman estimates and the Carleman type 
estimates for hyperbolic equations are subjected to so-called non-trapping 
conditions. Therefore both assumptions made in \cite{Kl} are critically 
important for the application of this method except of the  one dimensional case. 
In \cite{IY23}, the authors extended results of \cite{Kl} to the case of 
general second order hyperbolic equation. 
The main purpose of the current work is to remove the non-trapping assumptions
and prove the uniqueness without any geometric constraints on the 
observation subbondary $\gamma$.
\\

Henceforth we set 
$$
-A_1(x,D) = -A(x,D) + c_1(x), \quad -A_2(x,D) = -A(x,D) + c_2(x)
$$
with the domains $\DDD(A_1) = \DDD(A_2) = \DDD(A)$.
It is known that the spectrum 
$\sigma(A_k)$ of $A_k$, $k=1,2$, consists entirely of eigenvalues with 
finite multiplicities.

By changing $\www{u}:= e^{Mt}u$ with some constant $M$, it suffices to 
assume that there exists a constant $\kappa_2>1$ 
such that $(A_1u,u)_{L^2(\Omega)} \ge \kappa_2\Vert u\Vert^2_{L^2(\Omega)}$ 
for $u \in \DDD(A_1)$
and $(A_2v,v)_{L^2(\Omega)} \ge \kappa_2\Vert v\Vert^2_{L^2(\Omega)}$
for $v\in \DDD(A_2)$.

Then, setting
$$
\sigma(A_1)= \{\la_k\}_{k\in \N}, \quad \sigma(A_2) = \{ \mu_k\}_{k\in \N},
$$
we can number as
$$
1 < \la_1 < \la_2 < \cdots, \qquad 1<\mu_1 < \mu_2 < \cdots.
$$

Let $P_k$ be the projection for $\la_k$, $k \in \N$ which is defined by
$$
P_k = \frac{1}{2\pi\sqrt{-1}} \int_{\gamma(\la_k)}
(z-A_1)^{-1} dz, \quad
Q_k = \frac{1}{2\pi\sqrt{-1}} \int_{\gamma(\mu_k)}
(z-A_2)^{-1} dz,
$$
where $\gamma(\la_k)$ is a circle centered at $\la_k$ with sufficiently 
small radius such that the disc bounded by $\gamma(\la_k)$ does not
contain any points in $\sigma(A_1)\setminus \{\la_k\}$, and 
$\gamma(\mu_k)$ is a similar sufficiently small circle centered 
at $\mu_k$. 
Then $P_k:L^2(\OOO) \longrightarrow L^2(\OOO)$ is a bounded linear operator
to a finite dimensional space 
and $P_k^2 = P_k$ and $P_kP_{\ell} = 0$ for $k, \ell\in \N$ with 
$k \ne \ell$. Then 
$P_kL^2(\OOO) = \{ b\in \DDD(A_1);\, A_1b=\la_kb\}$, and we have
$a = \sum_{k=1}^{\infty} P_ka$ in $L^2(\OOO)$ for each 
$a \in L^2(\OOO)$ (e.g., Agmon \cite{Ag}, Kato \cite{Ka}).
Setting $m_k:= \mbox{dim}\, P_kL^2(\OOO)$, we have 
$m_k<\infty$, and we call $m_k$ the multiplicity of $\la_k$.
Similarly let $n_k$ and $Q_k$ be the multiplicity and
the eigenprojection for $\mu_k$. 
\\

Moreover we set $Q:= \OOO \times (0,T)$, and 
$$
H^{2,1}(Q):= \{ w\in L^2(Q);\, 
w, \, \ppp_iw,\, \ppp_i\ppp_jw,\, \ppp_tw \in L^2(Q)
\,\, \mbox{for $1\le i,j\le n$}\}.
$$
%Let $\Vert \cdot\Vert$ denote the norm in $L^2(\OOO)$, and we 
%specify the norm in a space $X$ by $\Vert \cdot\Vert_X$ if we 
%consider other space norms other than $L^2(\OOO)$.

Let 
$$
\Gamma=\{x\in \gamma;\, \vert a(x)\vert >0\}.
$$

We assume that 
$$
\Gamma\ne \emptyset. \eqno{(1.4)}
$$
For $a\in C(\ooo{\OOO})$, we set 
$$
\OOO_0 := \{ x\in \OOO;\, \vert a(x)\vert > 0\}.   
$$
%Let $\omega \subset\OOO_0$  be a union of open connected sets with the 
%following property: If $\tilde \omega$ is one of these sets  there exist an 
%open subset $\tilde \Gamma(\tilde \omega)\subset \overline{\tilde \omega}
%\cap \Gamma$ and 
%$$
%\vert a(x)\vert>0\quad \mbox{on}\quad  \tilde\Gamma(\tilde \omega).
%$$  

For $\Gamma$, we define 
$$
 \mbox{$\omega: = \{ x\in \OOO_0;\,$ there exist a point $x_*\in \Gamma$ 
and}
$$
$$
\mbox{a smooth curve $\ell \in C^\infty[0,1]$ such that  
$\ell(\xi) \in \OOO_0$ for $0<\xi\le 1$ and $\ell(0)=x_*$, $\ell(1)=x\}$}.
                                    \eqno{(1.5)}
$$
We remark that the definition (1.5) implies $\ell \setminus \{x_*\}
\subset \omega$.

In (1.5), replacing smooth curves by piecewise
smooth curves, we still have the same definition for $\omega$.

We note that $\omega$ is not necessarily a connected set.
However, if in addition we suppose that
$$
\mbox{$\Gamma$ is a connected subset of $\partial\Omega$},  \eqno{(1.6)}
$$
then one can verify that $\omega \subset \OOO$ is a domain, that is,
a connected open set.  Indeed, choosing $x, \www{x}\in \omega$ arbitrarily, 
we will show that we can find a piecewise smooth curve $L \subset \omega$ 
connecting $x$ and $\www{x}$ as follows.  
First we can choose smooth curves $\ell, \www{\ell}
\subset \OOO_0 \cup \Gamma$ and points $x_*, \www{x}_* \in \Gamma$
such that $\ell$ connects $x$ and $x_*$, $\www{\ell}$ connects $\www{x}$ and
$\www{x}_*$.  The definition implies that $\ell \setminus \{x_*\},
\www{\ell} \setminus \{ \www{x_*}\} \subset \omega$.
Since $\vert a\vert >0$ in $\Gamma$, we can 
find a smooth curve $\www{\gamma} \subset \OOO_0$ connecting $x_*$ and 
$\www{x}_*$.  Therefore, since $\www{\gamma} \subset \omega$, it follows that 
$x$ and $\www{x}$ can be connected by a piecewise smooth curve 
$L \subset \omega$ composed by $\ell, \www{\ell}, \www{\gamma}$,
which means that $\omega$ is a connected set.
Moreover, if $x\in \omega$, then we see that any point $\www{x} \in \OOO_0$ 
which is sufficiently close to $x$, can be connected to some point $\www{x}_*
\in \Gamma$ by some smooth curve in $\OOO_0$.  
Therefore, $\omega$ is a connected and 
open set, that is, $\omega$ is a domain.
$\blacksquare$

We can understand that $\omega$ is the maximal set such that 
all the points of $\omega$ is connected by a curve 
in $\OOO_0$ to $\Gamma$.
By (1.4), we note that
$\omega \ne \emptyset$. 
\\
{\bf Examples.}
\\
(i) Under condition (1.4), we have $\omega = \OOO$ if $\OOO_0 = \OOO$.
In general, if $\{ x\in \OOO;\, a(x) = 0\}$ has no interior points, then 
$\omega = \OOO_0$. 
\\
(ii) Assume that (1.4) and (1.6) hold true. Let subdomains $D_1, ..., D_m \subset \OOO$ satisfy 
$\ooo{D_1}, ..., \ooo{D_2} \subset \OOO$ and
$a=0$ on $\ooo{D_k}$ for $1\le k \le m$ and $\vert a\vert > 0$ in 
$\OOO \setminus \ooo{\bigcup_{k=1}^m D_k}$.  Then 
$\omega = \OOO \setminus \ooo{\bigcup_{k=1}^m D_k}$.
\\
(iii) Assume that (1.4) and (1.6) hold true. Let sudomains $D_1, D_2$ satisfy $\ooo{D_1} \subset D_2$, 
$\ooo{D_2} \subset \OOO$, $a=0$ in $\ooo{D_2 \setminus D_1}$ and
$\vert a\vert > 0$ in $D_1 \cup (\OOO \setminus \ooo{D_2})$.  Then 
$\omega = \OOO \setminus \ooo{D_2}$.  We note that $D_1$ is not included in 
$\omega$ although $\vert a\vert > 0$ in $D_1$.
\\

Now we state the main uniqueness result.
\\
{\bf Theorem 1.}
\\
{\it
Let $ a\in C(\ooo\Omega),$ $u,v \in H^{2,1}(Q)$ satisfy (1.2) and (1.3) 
respectively and $\ppp_tu, \ppp_tv \in H^{2,1}(Q)$ and let (1.4) hold true.
Assume
\\
{\bf Condition 1:} there exists a function $\theta \in C[1,\infty)$ satisfying
$$
\lim_{\eta\to\infty} \frac{\theta(\eta)}{\eta^{\frac{2}{3}}} = +\infty
$$
and
$$
\sumk e^{\theta(\la_k)} \Vert P_ka\Vert^2_{L^2(\Omega)} < \infty 
\quad \mbox{or}
\quad \sumk e^{\theta(\mu_k)} \Vert Q_ka\Vert^2_{L^2(\Omega)} < \infty.
                                         \eqno{(1.7)}
$$
Then,  
$$
u=v \quad \mbox{on $\gamma \times (0,T)$}
$$
implies $c_1=c_2$ on $\ooo{\omega}$.
}

As is seen by the proof, without the assumption (1.7), we can prove 
at least the coincidence of the eigenvalues of $A_1$ and $A_2$ of 
non-vanishing modes:
\\
{\bf Corollary.}
\\
{\it 
Let $u=v$ on $\gamma \times (0,T)$.  Then
$$
\{ \la_k;\, k\in \N, \, P_ka \ne 0 \quad \mbox{in $\OOO$}\}
= \{ \mu_k;\, k\in \N, \, Q_ka \ne 0 \quad \mbox{in $\OOO$}\}
$$
and if $P_ka \ne 0$ in $\OOO$ for $k\in \N$, then 
$$
P_ka = Q_ka \quad \mbox{on $\gamma$},
$$
after suitable re-numbering of $k$.
}

The corollary means that $u=v$ on $\gamma \times (0,T)$ implies 
that there exists $N_1 \in \N \cup \{\infty\}$ such that we can find
sequences $\{i_k\}_{1\le k\le N_1}, \, \{j_k\}_{1\le k\le N_1}
\subset \N$ satisfying 
$$
\left\{ \begin{array}{rl}
& \la_{i_k} = \mu_{j_k}, \quad P_{i_k}a \ne 0, \,\, 
Q_{j_k}a \ne 0 \quad \mbox{in $\OOO$},\quad
   P_{i_k}a = Q_{j_k}a = 0 \quad \mbox{on $\gamma$}
\quad \mbox{for $1\le k \le N_1$}, \\
& P_ia = 0 \quad\mbox{in $\OOO$ if $i\not\in \{ i_k\}_{1\le k\le N_1}$},
\quad
 Q_ja = 0 \quad\mbox{in $\OOO$ if $j\not\in \{ j_k\}_{1\le k\le N_1}$}.
\end{array}\right.
$$
We remark that even in the case $N_1=\infty$, we may have
$\{ i_k\}_{1\le k \le N_1} \subsetneqq \N$.
\\
{\bf Remark.}
In (1.7), consider a function $\theta(\eta) = \eta^p$.
Theorem 1 asserts the uniqueness 
if the initial value $a$ is smooth in the sense (1.7).
We emphasize that in (1.7), the critical exponent of $\la_k$ should be 
greater than $\frac{2}{3}$.  If we assume
the stronger condition $p=1$, that is,
$$
\sumk e^{\sigma \la_k} \Vert P_ka\Vert^2_{L^2(\Omega)} < \infty \quad \mbox{and}
\quad \sumk e^{\sigma \mu_k} \Vert Q_ka\Vert^2_{L^2(\Omega)} < \infty
                                         \eqno{(1.8)}
$$
with some constant $\sigma>0$, then the uniqueness is trivial 
because we can extend the solutions $u(\cdot,t)$ and $v(\cdot,t)$
to the time interval $(-\delta, 0)$ with small $\delta > 0$.
Indeed, since $\sumk \vert e^{\hhalf \sigma \la_k}\vert^2 \Vert P_ka\Vert^2_{L^2(\Omega)} 
< \infty$, we can verify that $u(\cdot,t) = \sumk e^{-\la_kt}P_ka$ in 
$L^2(\OOO)$ for $t > -\hhalf\sigma$.  Therefore we can extend 
$u(\cdot,t)$ to $\left(-\frac{\sigma}{2},\, 0\right)$ in $L^2(\OOO)$ and also
to $(-\delta, 0)$ with sufficiently small $\delta>0$.  The extension of
$v(\cdot,t)$ is similarly done.
Therefore, under (1.8), our inverse problem is reduced to the case where
the spatial data of $u,v$ are given at an intermediate time of the whole
time interval under consideration, 
which has been already solved in Bukhgeim and Klibanov
\cite{BK}, Imanuvilov and Yamamoto \cite{IY98}, Isakov \cite{Is}.

Condition corresponding to the case $p=\hhalf$ in (1.7) 
appears in the controllability of a parabolic equation.  
We know that a function $a(\cdot)$ in $\OOO$ satisfying the condition (1.7)
with $\theta(\eta) = \eta^{\hhalf}$ and $c_1\equiv 0$
belongs to the reachable set 
$$
\{ u(\cdot,0);\, b\in L^2(\OOO),\, h \in L^2(\ppp\OOO \times (-\tau,0)\},
$$
where $u$ is the solution to 
$$
\left\{ \begin{array}{rl}
& \ppp_tu = \Delta u \quad \mbox{in $\OOO \times (-\tau,0)$}, \\
& \ppp_{\nu}u = h \quad \mbox{on $\ppp\OOO \times (-\tau,0)$},\\
& u(\cdot,-\tau) = b \quad \mbox{in $\OOO$}
\end{array}\right.
$$
(Theorem 2.3 in Russell \cite{R}).
See also (1.9) stated below.
\\

The article is composed of four sections.  In Section 2, we show 
Carleman estimates for elliptic operator.  In Section 3, we prove the
uniqueness for our inverse problem first under a condition:
there exists a constant $\sigma_1>0$ such that 
$$
\sumk e^{\sigma_1 \la_k^{\hhalf}}\Vert P_ka\Vert^2_{L^2(\Omega)} 
+ \sumk e^{\sigma_1 \mu_k^{\hhalf}} \Vert Q_ka\Vert^2_{L^2(\Omega)} < \infty,
                                         \eqno{(1.9)}
$$
and next by proving that (1.7) yields (1.9), we complete
the proof of Theorem 1.
In Section 4, we prove a Carleman estimate used for deriving (1.9) from 
(1.7).
\section{Key Carleman estimate}

The proof of Theorem 1 relies essentially on the reduction of our
inverse parabolic problem to an inverse elliptic problem.
After the reduction, we prove the uniqueness by the method developed in 
\cite{BK} or, \cite{HIY}, \cite{IY98}, 
and so we need a relevant Carleman estimate 
for an elliptic equation.
For the statement of Carleman estimate, we introduce a weight function.

We arbitrarily fix $y \in \omega$. For $y$, we construct a non-empty
domain $\omega_y \subset \OOO$ satisfying
$$
\left\{ \begin{array}{rl}
&\mbox{(i)} \,\,y \in \omega_y, \quad \omega_y \subset \omega. \\
&\mbox{(ii)} \,\,\mbox{$\ppp\omega_y$ is of $C^{\infty}$-class.}\\
&\mbox{(iii)} \,\, \mbox{$\ppp\omega_y \cap \Gamma$ has interior points 
in the topology of $\ppp\OOO$.} \\
&\mbox{(iv)} \,\, \vert a(x)\vert > 0 \quad \mbox{for all 
  $x\in \ooo{\omega_y}$}.
\end{array}\right.
                                        \eqno{(2.1)}
$$
Indeed, since $y\in \omega$, by the definition of $\omega$, we can find
$y_*\in \Gamma$ and a smooth curve $\ell \in C^{\infty}[0,1]$ such that 
$\ell(1) = y$ and $\ell(0) = y_*$, $\ell(\xi) \in \omega$ for 
$0<\xi \le 1$.  Then as $\omega_y$, we can choose
a sufficiently thin neighborhood of the curve $\{\ell(\xi);\, 0<\xi\le 1\}$
which is included in $\omega$. 

For the proof of Theorem 1, we will show that if $y\in\omega$ then  $y\notin \mbox{supp}\, f.$  This of course implies that $f=0$ on $\omega$.
First we establish a Carleman estimate in $\omega_y \times (-\tau,\tau)$ with 
a constant $\tau > 0$.
 
We know that there exists a function
$d\in C^2(\ooo{\omega_y})$ such that 
$$
\vert \nabla d(x)\vert > 0 \quad \mbox{for $x\in \ooo{\omega_y}$}, \quad
d(x) > 0  \quad \mbox{for $x\in \omega_y$}, \quad
d(x) = 0 \quad \mbox{for $x\in \ppp\omega_y\setminus \Gamma$}.
                                                           \eqno{(2.2)}
$$
The existence of such $d$ is proved for example in Imanuvilov \cite{Im}.
See also Fursikov and Imanuvilov \cite{FI}.

For a constant $\tau>0$, we set 
$$
{\mathcal Q}_\tau:= \omega_y \times (-\tau, \, \tau),
$$
$\ppp_0 := \frac{\ppp}{\ppp t}$, and 
$$
\alpha(x,t) := e^{\la(d(x) - \beta t^2)}, 
\quad (x,t)\in {\mathcal Q}_{\tau}                           \eqno{(2.3)}
$$
with an arbitrarily chosen constant $\beta > 0$ and sufficiently large 
fixed $\la > 0$.
Then
\\
{\bf Lemma 2.1 (elliptic Carleman estimate).}
\\
{\it 
There exists a constant $s_0 > 0$ such that we can find a constant 
$C>0$ such that 
\begin{align*}
& \int_{{\mathcal Q}_\tau} \left\{ \frac{1}{s}\sum_{i,j=0}^n \vert \ppp_i\ppp_jw\vert^2
+ s\vert \ppp_tw\vert^2 + s\vert \nabla w\vert^2
+ s^3\vert w\vert^2\right) e^{2s\alpha} dxdt\\
\le& C\int_{{\mathcal Q}_\tau} \vert \ppp_t^2w - A_1w\vert^2 e^{2s\alpha} dxdt
\end{align*}
for all $s \ge s_0$ and $w\in H^2_0({\mathcal Q}_\tau)$.
}
\\

Here we recall that $-A_1w = \sumij \ppp_i(a_{ij}(x)\ppp_jw)
+ c_1(x)w$.   
The constants $s_0>0$ and $C>0$ can be chosen uniformly provided that 
$\Vert c_1\Vert_{L^{\infty}(\omega)} \le M$: arbitrarily fixed constant
$M>0$.

We note that Lemma 2.1 is a Carleman estimate for the elliptic operator 
$\ppp_t^2 - A_1w$.  Since $(\nabla \alpha, \, \ppp_t\alpha)
= (\nabla d, \, -2\beta t) \ne (0,0)$ on $\ooo{{\mathcal Q}_{\tau}}$ by (2.2), 
the proof of
the lemma relies directly on integration by parts and standard, similar for
example to the proof of Lemma 7.1 (p.186) in 
Bellassoued and Yamamoto \cite{BY}.
See also H\"ormander \cite{H}, Isakov \cite{Is}, where the estimation
of the second-order derivatives is not included but can be be 
derived by the a priori estimate for the elliptic boundary value problem.

For the proof of Theorem 1, we further need another Carleman estimate 
in $\OOO$ for an elliptic 
equation.  We can find $\rho\in C^2(\ooo{\OOO})$ such that 
$$
\rho(x) > 0 \quad \mbox{for $x \in\OOO$}, \quad
\vert \nabla \rho(x) \vert > 0 \quad \mbox{for $x \in\ooo{\OOO}$}, \quad
\ppp_{\nu_A}\rho(x) \le 0 \quad \mbox{for $x\in \ppp\OOO\setminus \gamma$}.
                                                      \eqno{(2.4)}
$$
The construction of $\rho$ can be found in Lemma 2.3 in \cite{IY98} for 
example.
Moreover, fixing a constant $\la>0$ large, we set 
$$
\psi(x):= e^{\la(\rho(x) - 2\Vert \rho\Vert_{C(\ooo{\OOO})})}, \quad 
x\in \OOO.
$$
Then
\\
{\bf Lemma 2.2.}
\\
{\it
There exist constants $s_0>0$ and $C>0$ such that 
$$
 \int_{\OOO} (s^3\vert g\vert^2 + s\vert \nabla g\vert^2)e^{2s\psi(x)} dx
\le C\int_{\OOO} \vert A_2g\vert^2 e^{2s\psi(x)} dx 
+ cs^3\int_{\gamma} (\vert g\vert^2 + \vert \nabla g\vert^2) e^{2s\psi} dS
$$
for all $s>s_0$ and $g \in H^2(\OOO)$ satisfying $\ppp_{\nu_A}g=0$ 
on $\ppp\OOO$.
}

We postpone the proof of Lemma 2.2 to Section 4.
\section{Proof of Theorem 1.1}

We divide the proof into four steps.
In Steps 1-3, we assume the condition (1.9) to prove the 
conclusion on uniqueness in Theorem 1.1.
\\
{\bf First Step.}
\\
We write $u(t):= u(\cdot,t)$ and $v(t):= v(\cdot,t)$ for $t>0$.
We recall that
$$
u(t) = \sumk e^{-\la_kt}P_ka, \quad
v(t) = \sumk e^{-\mu_kt}Q_ka  \quad \mbox{in $H^2(\OOO)$ for $t>0$.}
$$

We can choose subsets $\N_1, \M_1 \subset \N$ such that 
$$
\N_1:= \{k \in \N;\, P_ka \not\equiv 0 \quad \mbox{in $\OOO$}\}, \quad
\M_1:= \{k \in \N;\, Q_ka \not\equiv 0 \quad \mbox{in $\OOO$}\}.
                                           \eqno{(3.1)}
$$
We note that $\N_1 = \N$ or $\M_1 = \N$ may happen.

We can renumber the sets $\N_1$ and $\M_1$ as 
$$
\N_1 = \{1, ...., N_1\}, \quad \M_1=\{ 1, ...., M_1\},
$$
where $N_1 = \infty$ or $M_1 = \infty$ may occur.  
By $^{\ssss}\N_1$ we mean the cardinal number of the set $\N_1$.
Wote that 
$$
\la_1< \la_2 < \cdots < \la_{N_1} \quad \mbox{if $^{\ssss}\N_1 < \infty$} \quad
\la_1< \la_2 < \cdots \quad \mbox{if $^{\ssss}\N_1 = \infty$}
$$
and
$$
\mu_1< \mu_2 < \cdots < \mu_{M_1} \quad \mbox{if $^{\ssss}\M_1 < \infty$} \quad
\mu_1< \mu_2 < \cdots \quad \mbox{if $^{\ssss}\M_1 = \infty$}.
$$

Assuming that $u=v$ on $\gamma \times (0,T)$, by the time analyticity
of $u(t)$ and $v(t)$ for $t>0$ (e.g., Pazy \cite{Pa}), we obtain
$$
\sum_{k=1}^{N_1} e^{-\la_kt}P_ka = \sum_{k=1}^{M_1} e^{-\mu_kt}Q_ka 
\quad \mbox{on $\gamma \times (0,\infty)$}.    
                                                \eqno{(3.2)}
$$
We will prove that $\la_1 = \mu_1$.
Assume that $\la_1 < \mu_1$.  Then 
$$
P_1a + \sum_{k=2}^{N_1} e^{-(\la_k-\la_1)t}P_ka 
= \sum_{k=1}^{M_1} e^{-(\mu_k-\la_1)t}Q_ka 
\quad \mbox{on $\gamma \times (0,\infty)$}.
$$
Since $\la_k - \la_1 > 0$ for $2\le k \le N_1$ and 
$\mu_k - \la_1 > 0$ for $1\le k \le M_1$, letting $t\to \infty$,
we see that $P_1a=0$ on $\Gamma$.
Therefore
$$
\left\{ \begin{array}{rl}
& (A_1-\la_1)P_1a = 0 \quad \mbox{in $\OOO$}, \\
& P_1a\vert_{\Gamma} = 0, \quad \nuA P_1a\vert_{\ppp\OOO} = 0.
\end{array}\right.
$$
The unique continuation for the elliptic equation $A_1P_1a = \la_1P_1a$, (see e.g. \cite{H})
yields that
$$
P_1a = 0 \quad \mbox{in $\OOO$}.
$$
This is a contradiction by $1 \in \N_1$.
Thus  the inequality $\la_1 < \mu_1$ is impossible.
Similarly we can see that  the inequality $\la_1 > \mu_1$ is impossible. 
Therefore $\la_1 = \mu_1$ follows.

By (3.2) and $\la_1 = \mu_1$, we have
$$
P_1a - Q_1a 
= -\sum_{k=2}^{N_1} e^{-(\la_k-\la_1)t}P_ka 
+ \sum_{k=2}^{M_1} e^{-(\mu_k-\la_1)t}Q_ka 
\quad \mbox{on $\gamma \times (0,\infty)$}.
$$
Hence, by $\la_k-\la_1 > 0$ and $\mu_k - \la_1 = \mu_k - \mu_1 > 0$ for 
all $k \ge 2$, letting $s \to \infty$ we obtain $P_1a = Q_1a$ on $\gamma$.

In view of (3.2), we obtain
$$
\sum_{k=2}^{N_1} e^{-\la_kt}P_ka 
= \sum_{k=2}^{M_1} e^{-\mu_kt}Q_ka 
\quad \mbox{on $\gamma \times (0,\infty)$}. 
$$
Repeating the same argument as much as possible, we reach 
$$
N_1 = M_1, \quad \la_k = \mu_k, \quad
P_ka = Q_ka \quad \mbox{on $\gamma \times (0,\infty)$ for
$1\le k \le N_1$}.                              \eqno{(3.3)}
$$
\\
{\bf Second Step.}

We consider two initial boundary value problems for elliptic equations:
$$
\left\{ \begin{array}{rl}
& \ppp_t^2w_1 - A_1w_1 = 0 \quad \mbox{in $\OOO\times (0,\tau)$}, \\
& \nuA w_1 = 0 \quad \mbox{on $\ppp\OOO \times (0,\tau)$}, \\
& w_1(x,0) = a(x), \quad \ppp_tw_1(x,0) = 0, \quad x\in \OOO
\end{array}\right.
                                 \eqno{(3.4)}
$$
and
$$
\left\{ \begin{array}{rl}
& \ppp_t^2w_2 - A_2w_2 = 0 \quad \mbox{in $\OOO\times (0,\tau)$}, \\
& \nuA w_2 = 0 \quad \mbox{on $\ppp\OOO \times (0,\tau)$}, \\
& w_2(x,0) = a(x), \quad \ppp_tw_2(x,0) = 0, \quad x\in \OOO.
\end{array}\right.
                                \eqno{(3.5)}
$$
Since we have the spectral representations by (1.9), we can obtain 
$$
e^{-tA_1^{\hhalf}}a = \sumk e^{-\la_k^{\hhalf}t}P_ka, \quad
e^{tA_1^{\hhalf}}a = \sumk e^{\la_k^{\hhalf}t}P_ka \quad 
\mbox{in $L^2(\OOO)$ for $t>0$}
$$
and similar representations hold for $e^{\pm tA_2^{\hhalf}}a$.
 
Then by the assumption (1.9) on $a$, we see that 
$$
w_1(t) = \frac{1}{2}(e^{-tA_1^{\hhalf}}a + e^{tA_1^{\hhalf}}a), \quad
w_2(t) = \frac{1}{2}(e^{-tA_2^{\hhalf}}a + e^{tA_2^{\hhalf}}a) 
$$
in $H^2(\Omega\times (0,\tau))$ for $t \in (0,\tau)$, satisfy (3.4) and (3.5)
respectively if $\tau>0$ is chosen sufficiently small.

In view of (3.3) and the definition (3.1) of $\N_1$, 
the spectral representations imply
\begin{align*}
& w_1(x,t) = \hhalf \sum_{k=1}^{N_1} (e^{-\la_k^{\hhalf}t}P_ka
+ e^{\la_k^{\hhalf}t}P_ka)
+ \hhalf \sum_{k\in \N \setminus \{1, ..., N_1\}}(e^{-\la_k^{\hhalf}t}P_ka
+ e^{\la_k^{\hhalf}t}P_ka)\\
=& \hhalf \sum_{k=1}^{N_1} (e^{-\la_k^{\hhalf}t}P_ka
+ e^{\la_k^{\hhalf}t}P_ka) \quad \mbox{in $\OOO\times (0,\tau)$,}
\end{align*}
and 
$$
w_2(x,t) = \hhalf \sum_{k=1}^{N_1} (e^{-\la_k^{\hhalf}t}Q_ka
+ e^{\la_k^{\hhalf}t}Q_ka) \quad \mbox{in $\OOO\times (0,\tau)$}.
$$
Therefore (3.3) yields
$$
w_1 = w_2 \quad \mbox{on $\gamma \times (0,\tau)$.}       \eqno{(3.6)}
$$

Now we reduce our inverse problem for the parabolic equations to the one
for elliptic equations for (3.4) and (3.5).
This is the essence of the proof.
\\
{\bf Third Step.}
\\
By (1.9), we can readily verify further regularity $\ppp_tw_1, \ppp_tw_2
\in H^2(0,\tau;H^2(\OOO))$.
Setting $y:= w_1 - w_2$ and $R:=w_2$ in $\OOO\times (0,\tau)$ 
and $f:= c_2-c_1$ in $\OOO$, by (3.4) -(3.6) we have
$$
\left\{ \begin{array}{rl}
& \ppp_t^2y - A_1y = f(x)R(x,t) \quad \mbox{in $\OOO\times (0,\tau)$}, \\
& \nuA y = 0 \quad \mbox{on $\ppp\OOO \times (0,\tau)$}, \\
& y=0 \quad \mbox{on $\gamma \times (0,\tau)$}, \\
& y(x,0) = \ppp_ty(x,0) = 0, \quad x\in \OOO.
\end{array}\right.
                                                   \eqno{(3.7)}
$$

Now we will prove that for any pair $(y,f)$ solving problem (3.7) we have $y\notin \mbox{supp}, f.$ Since $ y$ was chosen as an arbitrary point from $\omega$ this implies
$$
f=0\quad \mbox{in}\quad \omega.
$$  The argument relies on \cite{IY98}.

We set
$$
\www{y}(x,t) = 
\left\{ \begin{array}{rl}
& y(x,t), \quad 0<t<\tau, \\
& y(x,-t), \quad -\tau<t<0
\end{array}\right.
$$
and
$$
\www{R}(x,t) = 
\left\{ \begin{array}{rl}
& R(x,t), \quad 0<t<\tau, \\
& R(x,-t), \quad -\tau<t<0.
\end{array}\right.
$$
Then we see that $\www{y} \in H^3(-\tau,\tau;H^2(\omega_y))$ by
$y(\cdot,0) = \ppp_ty(\cdot,0) = 0$ in $\omega_y$, and 
$$
\left\{ \begin{array}{rl}
& \ppp_t^2\www{y} - A_1\www{y} = f(x)\www{R}(x,t) 
\quad \mbox{in ${\mathcal Q}_{\tau}:= \omega_y\times (-\tau,\tau)$}, \\
& \www{y}= \nuA \www{y} = 0 \quad \mbox{on $(\ppp\omega_y \cap \Gamma) 
\times (-\tau,\tau)$}, \\
& \www{y}(x,0) = \ppp_t\www{y}(x,0) = 0, \quad x\in \omega_y.
\end{array}\right.
$$
Setting $\www{z}:= \ppp_t\www{y}$, we have $\www{z} \in H^2(Q_{\tau})$ and
$$
\left\{ \begin{array}{rl}
& \ppp_t^2\www{z} - A_1\www{z} = f(x)\ppp_t\www{R}(x,t) 
\quad \mbox{in ${\mathcal Q}_{\tau}$}, \\
& \www{z} = \nuA \www{z} = 0 \quad \mbox{on $(\ppp\omega_y \cap \Gamma)
\times (-\tau,\tau)$}, \\
& \www{z}(x,0) = 0, \quad x\in \omega_y
\end{array}\right.
                                     \eqno{(3.8)}
$$
and
$$
\ppp_t\www{z}(x,0) = f(x)a(x), \quad x\in \omega_y.     \eqno{(3.9)}
$$
In the weight function $d(x) - \beta t^2$ in Lemma 2.1, for $\tau>0$ we choose 
$\beta > 0$ sufficiently large,
so that 
$$
\Vert d\Vert_{C(\ooo{\omega_y})} - \beta \tau^2 < 0,   \eqno{(3.10)}
$$
We choose constants $\delta_1, \delta_2 > 0$ such that
$$
\Vert d\Vert_{C(\ooo{\omega_y})} - \beta \tau^2 < 0 < \delta_1 < \delta_2 
   \quad \mbox{and}\quad d(y)>\delta_2                           \eqno{(3.11)}
$$
and we define $\chi \in C^{\infty}(\ooo{{\mathcal Q}_{\tau}})$ satisfying
$$
\chi(x,t) = 
\left\{ \begin{array}{rl}
& 1, \quad d(x) - \beta t^2 > \delta_2, \\
& 0, \quad d(x) - \beta t^2 < \delta_1. 
\end{array}\right.
                                              \eqno{(3.12)}
$$
In particularl, $d(y) > \delta_2$ implies 
$$
\chi(y,0)=1.
$$
%The parameter $\ep>0$ will be used later for approximating $\omega$ 
%by $\{x\in \omega;\, d(x) > \ep\}$.
 
Setting $z:= \chi \www{z}$ on $\ooo{{\mathcal Q}_{\tau}}$, we see that 
$$
z = \nuA z = 0 \quad \mbox{on $\ppp\omega_y \times (-\tau, \tau)$}
                                                  \eqno{(3.13)}
$$
and
$$
z = \ppp_t z = 0 \quad \mbox{on $\omega_y \times \{ \pm \tau\}$}.
                                                \eqno{(3.14)}
$$
Indeed, $(x,t) \in (\ppp\omega_y \setminus \Gamma) \times (-\tau, \tau)$ 
implies
$d(x) - \beta t^2 = -\beta t^2 \le 0 < \delta_1$ by (2.2), and so 
the definition (3.12) of $\chi$ yields that 
$\chi(x,t) = 0$ in a neighborhood of such $(x,t)$.
For $(x,t) \in (\ppp\omega_y \cap \Gamma) \times (-\tau,\tau)$, 
by (3.8) we see that 
$z(x,t) = \nuA z(x,t) =  0$, which verifies (3.13).
Moreover, on $\omega_y\times \{ \pm \tau\}$, by (3.11) we have 
$$
d(x) - \beta t^2 \le \Vert d\Vert_{C(\ooo{\omega_y})} - \beta\tau^2
< 0 < \delta_1,
$$
so that $\chi(x,t) = 0$ in a neighborhood of such $(x,t)$.  
Thus (3.14) has been verified.
$\blacksquare$

Consequently we prove that $z\in H^2_0(Q_{\tau})$.
Moreover, we can readily obtain
$$
\left\{ \begin{array}{rl}
& \ppp_t^2z - A_1z = \chi(\ppp_t\www {R})f + R_0(x,t), \quad 
(x,t) \in {\mathcal Q}_{\tau}, \\
& z = \vert \nabla z\vert = 0 \quad \mbox{on $\ppp {\mathcal Q}_{\tau}$},
\end{array}\right.
                                                   \eqno{(3.15)}
$$
where $R_0$ is a linear combination of $\nabla \www{z}$, $\ppp_t\www{z}$,
whose coefficients are linear combinations of $\nabla\chi$ and 
$\ppp_t\chi$.  Therefore (3.12) implies 
$$
R_0(x,t) \ne 0 \quad \mbox{only if $\delta_1 \le d(x) - \beta t^2
\le \delta_2$}.                              \eqno{(3.16)}
$$

Therefore we can apply Lemma 2.1 to (3.15) with (3.16):
$$
\int_{{\mathcal Q}_{\tau}} \left( \frac{1}{s}\vert \ppp_t^2z\vert^2
+ s\vert \ppp_tz\vert^2 \right) e^{2s\alpha} dxdt 
                                                          \eqno{(3.17)}
$$
\begin{align*}
\le & C\int_{{\mathcal Q}_{\tau}} \vert \chi(\ppp_t\www{R})f \vert^2e^{2s\alpha} dxdt
+ C\int_{{\mathcal Q}_{\tau}} \vert R_0(x,t)\vert^2 e^{2s\alpha} dxdt \\
\le & C\int_{{\mathcal Q}_{\tau}}\chi^2 \vert f\vert^2 e^{2s\alpha} dxdt 
+ Ce^{2se^{\la\delta_2}}
\end{align*}
for all large $s>0$.

On the other hand, since $\ppp_tz(\cdot,-\tau) = 0$ in $\omega_y$ by 
(3.14), we have
\begin{align*}
& \int_{\omega_y} \vert \ppp_tz(x,0)\vert^2 e^{2s\alpha(x,0)} dx 
= \int^0_{-\tau} \ppp_t\left( \int_{\omega_y} \vert \ppp_tz(x,t)\vert^2
e^{2s\alpha(x,t)} dx \right) dt\\
=& \int^0_{-\tau} \int_{\omega_y} \{ 2(\ppp_tz)(x,t)\ppp_t^2z(z,t)
+ \vert \ppp_tz \vert^2 2s(\ppp_t\alpha))\} e^{2s\alpha(x,t)} dxdt.
\end{align*}
Since 
$$
\vert (\ppp_tz)(\ppp_t^2z)\vert \le \frac{1}{2}
\left( s\vert \ppp_tz\vert^2 + \frac{1}{s}\vert \ppp_t^2z\vert^2
\right) \quad \mbox{in $\mathcal{Q}_\tau$},
$$
in terms of (3.17) we obtain
$$
 \int_{\omega_y} \vert \ppp_tz(x,0)\vert^2 e^{2s\alpha(x,0)} dx 
\le C\int_{{\mathcal Q}_{\tau}} \left( \frac{1}{s} \vert \ppp_t^2z(x,t)\vert^2 
+ s\vert \ppp_tz\vert^2 \right) e^{2s\alpha} dxdt 
                                                  \eqno{(3.18)}
$$
$$
\le C\int_{{\mathcal Q}_{\tau}}\vert \chi\vert^2\vert f\vert^2 e^{2s\alpha} 
dxdt + Ce^{2se^{\la\delta_2}}
$$
for all large $s>0$.

Moreover, we have
\begin{align*}
& \ppp_tz(x,0) = \ppp_t(\chi\www{z})(x,0)
= (\ppp_t\chi)(x,0)\www{z}(x,0) + \chi(x,0)\ppp_t\www{z}(x,0)\\
=& \chi(x,0)f(x)a(x)
\end{align*}
by (3.9) and $\www{z}(x,0) = 0$ for $x \in \omega_y$ in (3.8).
Therefore, in terms of (2.1)-(iv), we obtain
$$
\vert \ppp_tz(x,0)\vert \ge C\vert \chi(x,0)f(x)\vert, \quad 
x\in \ooo{\omega_y}.
$$
Consequently (3.18) implies
$$
\int_{\omega_y} \vert\chi(x,0)\vert^2 \vert f(x) \vert^2 e^{2s\alpha(x,0)} dx 
\le C\int_{\mathcal{Q}_{\tau}} \vert \chi f\vert^2 e^{2s\alpha} dxdt 
+ Ce^{2se^{\la\delta_2}}             \eqno{(3.19)}   
$$
for all large $s>0$.

Moreover, we see
\begin{align*}
&\int_{{\mathcal Q}_{\tau}} \vert \chi f\vert^2 e^{2s\alpha} dxdt 
= \int^{\tau}_{-\tau} \int_{\omega_y} \vert\chi(x,0)f(x)\vert^2 e^{2s\alpha} 
dxdt\\
= & \int_{\omega_y} \vert \chi(x,0) f(x)\vert^2 e^{2s\alpha(x,0)} 
\left( \int^{\tau}_{-\tau} e^{2s(\alpha(x,t) - \alpha(x,0))} dt\right)dx.
\end{align*}
Since 
$$
\int^{\tau}_{-\tau} e^{2s(\alpha(x,t) - \alpha(x,0))} dt
= \int^{\tau}_{-\tau} e^{2se^{\la d(x)}(e^{-\la\beta t^2} - 1)} dt
\le \int^{\tau}_{-\tau} e^{Cs(e^{-\la\beta t^2} - 1)} dt
= o(1)
$$
as $s \to \infty$ by the Lebesgue convergence theorem.
Hence, (3.19) yields
$$
 \int_{\omega_y} \vert \chi(x,0)f(x)\vert^2 e^{2s\alpha(x,0)} dx
\le o(1)\int_{\omega_y} \vert \chi(x,0)f(x)\vert^2 e^{2s\alpha(x,0)} dx
+ Ce^{2se^{\la\delta_2}},
$$
and we can absorb the first term on the right-hand side into the 
left-hand side to reach
$$
\int_{\omega_y} \vert \chi(x,0)f(x)\vert^2 e^{2s\alpha(x,0)} dx
\le Ce^{2se^{\la\delta_2}}             \eqno{(3.20)}
$$
for all large $s>0$. 

Henceforth we set $B(y,\ep):= \{ x;\, \vert x-y\vert < \ep\}$.
Then, we can choose $\delta_3 > \delta_2$ and
a sufficiently small $\ep > 0$ such that
$B(y,\ep) \subset \omega_y$ and $d(x) \ge \delta_3$ for all 
$x\in B(y,\ep)$.  This is possible, because $d(y) > \delta_2$ in 
(3.11) and $\omega_y$ is an open set including $y$.

We shrink the integration region of the left-hand side of (3.20)
to $B(y,\ep)$ and obtain
$$
\int_{B(y,\ep)} \vert \chi(x,0)f(x)\vert^2 e^{2s\alpha(x,0)} dx
\le Ce^{2se^{\la\delta_2}}
$$
for all large $s>0$. 
Since $d(x) \ge \delta_3 > \delta_2$ for $x\in B(y,\ep)$, the condition (3.12)
yields $\chi(x,0) = 1$ and $\alpha(x,0) = e^{\la d(x)}
\ge e^{\la\delta_3}$ for all $x\in B(y,\ep)$.  Therefore,
$$
\left( \int_{B(y,\ep)} \vert f(x)\vert^2 dx\right) e^{2se^{\la\delta_3}}
\le Ce^{2se^{\la\delta_2}},
$$
that is, $\Vert f\Vert^2_{L^2(B(y,\ep))} \le e^{-2s(e^{\la\delta_3}
- e^{\la\delta_2})}$ for all large $s>0$. 

In terms of $\delta_3 > \delta_2$, letting $s\to \infty$, 
we see that the right-hand side
tends to $0$, and so $f=0$ in $B(y,\ep)$.  Since $y$ is arbitrarily chosen, 
we reach $f=c_2-c_1 = 0$ in $\omega$.
Thus the conclusion of Theorem 1.1 is proved under condition (1.9).
$\blacksquare$
\\
{\bf Fourth Step.}
\\
We will complete the proof of Theorem 1.1 by demonstrating
that (1.7) implies (1.9).
Without loss of generality, we can assume
$$
\sumk e^{\theta(\la_k)}\Vert P_ka\Vert^2_{L^2(\Omega)} < \infty.
                                     \eqno{(3.21)}
$$
It suffices to prove that there exists a constant $\sigma_1>0$ such that 
$$
\sum_{k=1}^{\infty} e^{\sigma_1\la_k^{\hhalf}}\Vert Q_ka\Vert^2_{L^2(\Omega)} < \infty,
                                                     \eqno{(3.22)}
$$
with the assumption that the set of $k\in \N$ such that 
$Q_ka\ne 0$ in $\OOO$ is infinite.
For simplicity, we can consider the case where 
$P_ka \ne 0$ in $\OOO$ for all $k\in \N$.  We can argue similarly 
in the rest cases.
Then, by Corollary which was already proved in First Step, we choose 
a subset $\M_1 \subset \N$ such that 
$$
\{ \la_i\}_{i\in\N} = \{ \mu_j\}_{j\in \M_1}, \quad
Q_ja = 0 \quad \mbox{in $\OOO$ for $j\in \N \setminus \M_1$}.
$$
Now it suffices to prove (3.22) in the case where $Q_ka\ne 0$ in 
$\OOO$ for all $k\in \N$.

After re-numbering, we can obtain
$$
\la_k = \mu_k, \quad P_ka = Q_ka \quad \mbox{on $\gamma$ for all
$k\in \N$}.                              \eqno{(3.23)}
$$
The trace theorem and the a priori estimate for an elliptic operator
yields
$$
\Vert P_ka\Vert_{H^1(\Gamma)} \le C\Vert P_ka\Vert_{H^2(\OOO)}
\le C(\Vert A_1P_ka\Vert_{L^2(\Omega)} + \Vert P_ka\Vert_{L^2(\Omega)})
= C(\la_k+1)\Vert P_ka\Vert_{L^2(\Omega)}.           \eqno{(3.24)}
$$
Here and henceforth $C>0$ denotes generic constants which are
independent of $s>0$ and $k\in \N$.

Since $A_2Q_k = \la_kQ_k$, by (3.23) we apply Lemma 2.2 to have
$$
s^3\int_{\OOO} \vert Q_ka\vert^2 e^{2s\psi} dx 
\le C\int_{\OOO} \la_k^2\vert Q_ka\vert^2 e^{2s\psi} dx
+ Cs^3\int_{\gamma} (\vert Q_ka\vert^2 + \vert \nabla (Q_ka)\vert^2)
e^{2s\psi} dx                               \eqno{(3.25)}
$$
$$
\le C\int_{\OOO} \la_k^2\vert Q_ka\vert^2 e^{2s\psi} dx
+ Cs^3e^{2sM}\Vert P_ka\Vert^2_{H^1(\gamma)}
$$
for all large $s>0$.  Here we set 
$M:= \max_{x\in \ooo{\Gamma}} \psi(x)$.

We choose $s>0$ sufficiently large and set 
$s_k:= s^*\la_k^{\frac{2}{3}}$ for $k\in \N$.  Then, using (3.24),
we obtain
\begin{align*}
& ({s^*}^3\la_k^2 - C\la_k^2)\int_{\OOO} \vert Q_ka\vert^2
e^{2s_k\psi} dx
\le C{s^*}^3\la_k^2e^{2s_kM}\Vert P_ka\Vert^2_{H^1(\OOO)}\\
\le& C{s^*}^3\la_k^2e^{2s_kM}(\la_k+1)^2\Vert P_ka\Vert^2_{L^2(\Omega)}.
\end{align*}
Since $\psi \ge 0$ in $\OOO$ and we can take $s^*>0$ sufficiently large, we see
$$
 {s^*}^3\la_k^2 \Vert Q_ka\Vert^2_{L^2(\Omega)}
\le C{s^*}^3\la_k^2\la_k^2 e^{2s_kM}\Vert P_ka\Vert^2_{L^2(\Omega)},
$$
that is,
$$
\Vert Q_ka\Vert^2_{L^2(\Omega)} \le C\la_k^2 e^{C_1\la_k^{\frac{2}{3}}}
\Vert P_ka\Vert^2_{L^2(\Omega)},
$$
where we set $C_1:= 2s^*M$.  Here we note that $s^*$ and $M$, and so
the constant $C_1$ are independent of $k\in \N$.

Therefore, since we can find a constant $C_2>0$ such that 
$\eta^2 e^{C_1\eta^{\frac{2}{3}} + \sigma_1\eta^{\hhalf}}
\le C_2e^{C_2\eta^{\frac{2}{3}}}$ for all $\eta \ge 0$, 
we see
$$
\sumk e^{\sigma_1\la_k^{\hhalf}}\Vert Q_ka\Vert^2_{L^2(\Omega)}
\le C\sumk \la_k^2e^{C_1\la_k^{\frac{2}{3}}+\sigma_1\la_k^{\hhalf}}
\Vert P_ka\Vert^2_{L^2(\Omega)}
\le C_2\sumk e^{C_2\la_k^{\frac{2}{3}}}\Vert P_ka\Vert^2_{L^2(\Omega)}.
$$
Moreover, $\lim_{k\to\infty} \frac{\theta(\la_k)}{\la_k^{\frac{2}{3}}}
= \infty$ yields that for the constant $C_2>0$ we can choose $N\in \N$ 
such that $C_2\la_k^{\frac{2}{3}} \le \theta(\la_k)$ for $k \ge N$.  
Consequently, 
$$
\sum_{k=N}^{\infty} e^{\sigma_1\la_k^{\hhalf}}\Vert Q_ka\Vert^2_{L^2(\Omega)}
\le C\sum_{k=N}^{\infty} e^{\theta(\la_k)}\Vert P_ka\Vert^2_{L^2(\Omega)}
< \infty,
$$
and so
$$
\sum_{k=1}^{\infty} e^{\sigma_1\la_k^{\hhalf}}\Vert Q_ka\Vert^2_{L^2(\Omega)}
< \infty.
$$
Thus (3.21) completes the proof of Theorem 1.1.
$\blacksquare$
\section{Appendix: Proof of Lemma 2.2}

We can prove the lemma by integration by parts similarly to 
Lemma 7.1 (p.186) in Bellassoued and Yamamoto \cite{BY}) for example, but
here we derive from a Carleman estimate for the parabolic equation
by Imanuvilov \cite{Im}.

We set $Q:= \OOO\times (0,T)$.
We choose $\ell \in C^{\infty}[0,T]$ such that 
$$
\left\{ \begin{array}{rl}
& \ell(t) = 1 \quad \mbox{for $\frac{T}{4}\le t\le \frac{3}{4}T$},\\
& \ell(0) = \ell(T) = 0,\\
& \mbox{$\ell$ is strictly increasing on $\left[ 0, \, \frac{T}{4}\right]$
and strictly decreasing on $\left[ \frac{T}{4}, \,T\right]$}.
\end{array}\right.
                                \eqno{(4.1)}
$$
In particular, $\ell(t) \le 1$ for $0\le t\le T$.  Choosing $\la>0$ 
sufficiently large, we set
$$
\alpha(x,t) := \frac{e^{\la\rho(x)} - e^{2\la\Vert \rho\Vert_{C(\ooo{\OOO})}}}
{\ell(t)}, \quad
\va(x,t) := \frac{e^{\la\rho(x)}}{\ell(t)}, \quad (x,t)\in \OOO\times (0,T).
$$
Then we know
\\
{\bf Lemma 4.1}
\\
{\it
There exist constants $s_0>0$ and $C>0$ such that 
\begin{align*}
& \int_Q (s\va\vert \nabla U\vert^2 + s^3\va^3\vert U\vert^2)
e^{2s\alpha} dxdt \\
\le& C\int_Q \vert \ppp_tU - A_2U\vert^2 e^{2s\alpha} dxdt
+ C\int_{\gamma} (\vert \ppp_tU\vert^2 + s\va\vert \nabla U\vert^2
+ s^3\va^3\vert U\vert^2) e^{2\alpha} dSdt
\end{align*}
for all $s \ge s_0$ and $U\in H^{2,1}(Q)$ satisfying 
$\nuA U = 0$ on $\ppp\OOO \times (0,T)$.
}

The proof is found in Chae, Imanuvilov and Kim \cite{CIK}.

We apply Lemma 4.1 to $g(x)$ satisfying $\nuA g = 0$ on $\ppp\OOO$ to
obtain
$$
\int_Q (s\va(x,t)\vert \nabla g(x)\vert^2 + s^3\va^3(x,t)\vert g(x)\vert^2)
e^{2s\alpha(x,t)} dxdt 
                                       \eqno{(4.2)}
$$
$$
\le C\int_Q \vert A_2g\vert^2 e^{2s\alpha(x,t)} dxdt
+ C\int^T_0 \int_{\gamma} (s\va(x,t)\vert \nabla g(x)\vert^2
+ s^3\va^3\vert g(x)\vert^2) e^{2\alpha} dSdt
$$
for all $s \ge s_0$.
Moreover, in terms of (4.1) and $e^{\la\rho(x)} 
= e^{2\la\Vert \rho\Vert_{C(\ooo{\OOO})}}\psi(x)$ for $x\in \OOO$,
we have 
$$
\int_Q (s\va(x,t)\vert \nabla g(x)\vert^2 + s^3\va^3(x,t)\vert g(x)\vert^2)
e^{2s\alpha(x,t)} dxdt
                                        \eqno{(4.3)}
$$
%\begin{align*}
$$
\ge \int^{\frac{3}{4}T}_{\frac{T}{4}} \int_{\OOO}
 (se^{\la\rho(x)}\vert \nabla g(x)\vert^2 
+ s^3e^{3\la\rho(x)}\vert g(x)\vert^2) 
\exp( 2s(e^{\la\rho(x)} - e^{2\la\Vert \rho\Vert_{C(\ooo\OOO)}}) ) dxdt
$$
%\\
$$
\ge C\frac{T}{2}
\int_{\OOO}  (s\vert \nabla g(x)\vert^2 + s^3\vert g(x)\vert^2) 
e^{2s(e^{2\la\Vert\rho\Vert_{C(\ooo\OOO)}}\psi(x))} dx 
e^{-2se^{2\la\Vert \rho\Vert_{C(\ooo\OOO)}}}
$$
$$
\ge C\frac{T}{2}
\int_{\OOO}  (s\vert \nabla g(x)\vert^2 + s^3\vert g(x)\vert^2) 
e^{2s\psi(x)} dx 
e^{-2se^{2\la\Vert \rho\Vert_{C(\ooo\OOO)} }}.
$$
%\end{align*}
Here $C>0$ depends on $\la$ but not on $s>0$.
By $e^{2s\alpha(x,t)} \le 1$ in $Q$, (4.2) and (4.3), we obtain
$$
C\frac{T}{2}\int_{\OOO}  (s\vert \nabla g(x)\vert^2 + s^3\vert g(x)\vert^2) 
e^{2s\psi(x)} dx e^{-2se^{2\la\Vert \rho\Vert_{C(\ooo\OOO)}}}
                                        \eqno{(4.4)}
$$
$$
\le C\int_Q \vert A_2g\vert^2 e^{2s\alpha(x,t)} dxdt
+ C\int^T_0 \int_{\gamma} (s\va(x,t)\vert \nabla g(x)\vert^2
+ s^3\va^3(x,t)\vert g(x)\vert^2) 
e^{2s\alpha(x,t)} dSdt.
$$
Since $\sup_{(x,t)\in \gamma \times (0,T)} 
\vert (s\va)^ke^{2s\alpha(x,t)}\vert 
< \infty$ for $k=1,3$ and 
$$
e^{2s\psi(x)} = e^{2se^{\la(\rho(x) - 2\Vert \rho\Vert_{C(\ooo{\OOO})})}}
\ge e^{2se^{-\la\Vert \rho\Vert_{C(\ooo\OOO)}}},
$$
we can find a constant $C_1 = C_1(\la) > 0$ such that 
$$
(s\va(x,t))^k e^{2s\alpha(x,t)} \le C_1e^{2s\psi(x)}, \quad
(x,t) \in \gamma \times (0,T).
$$
Therefore
\begin{align*}
& \int_{\OOO} (s\vert \nabla g\vert^2 + s^3\vert g\vert^2)
e^{2s\psi(x)} dx \\
\le &Ce^{2se^{2\la\Vert \rho\Vert_{C(\ooo{\OOO})}}}
\left( \int_{\OOO} \vert A_2g\vert^2 e^{2s\psi} dx 
+ \int_{\gamma} (s\vert \nabla g\vert^2 + s^3\vert g(x)\vert^2) 
e^{2s\psi} dS\right).
\end{align*}
Substituting (4.3) and (4.4) into (4.2), we complete the proof of 
Lemma 2.2.
$\blacksquare$

{\bf Acknowledgements.}
The work was supported by Grant-in-Aid for Scientific Research (A) 20H00117 
of Japan Society for the Promotion of Science.
\end{document}